\documentclass[10pt,a4paper]{article}
\usepackage{amsmath}
\usepackage{amssymb}
\usepackage{theorem}
\usepackage[utf8]{inputenc}
\theoremheaderfont{\scshape}
\newtheorem{thm}{Theorem}[section]

\newtheorem{lem}[thm]{Lemma}

\theorembodyfont{\normalfont}
\newtheorem{rem}[thm]{Remark}
\newcommand{\bb}[1]{\mathbb{#1}}
\newcommand{\ds}{\displaystyle}

\newcommand{\g}[1]{\mathbf{#1}}
\newcommand{\abs}[1]{\left\vert #1\right\vert}

\newcommand{\Lim}[2]{\lim\limits_{#1\rightarrow #2}}

\newcommand{\interventier}[2]{\ensuremath{[\![#1,#2]\!]}}

\title{Sharp asymptotics for the free energy of 1+1 dimensional directed polymers in an  infinitely divisible environment}

\date{}

\begin{document}

\maketitle

\begin{abstract}
We give sharp estimate for the free energy of directed polymers in random environment in dimension 1+1. This estimate was known for a Gaussian environment, we extend it to the case where the law of the environment is infinitely divisible.
\end{abstract}



\section{Introduction}

We refer to \cite{CometsShigaYoshidaReview} for a review of directed polymers in
random environment.
Let $S=(S_n)_{n\in\bb {N}}$ be the simple random walk on 
$\bb Z^d$ starting at $0$, defined on the probability space
$(\Sigma,\cal E,\textbf{P})$. Let $\eta=(\eta(n,x))_{(n,x)\in\bb N\times \bb
Z^d}$ be a sequence of real-valued i.i.d. random
variables defined on another probability space $(\Omega,\cal F, \bb P)$.
The expectation of a function $f$ with respect to the probability measure 
$\textbf{P}$ (respectively $\bb P$) will be denoted by $\textbf{E}f=\int_\Sigma f d\textbf{P}$
(respectively $\bb E f=\int_\Omega f d\bb P$).
The path $(i,S_i)_{1\leq i\leq n}$ represents the directed polymer of size $n$
in dimension $1+d$, and $\eta$ is the random environment. 
For $\beta$ strictly positive, featuring the inverse of the temperature, we define the random polymer probability measure
$\g P_{n,\beta}$ on the path space $(\Sigma,\cal E)$ by its density with respect to $\g P$
\begin{equation}
\frac{d\g P_{n,\beta}}{d\textbf{P}}(S)=\ds\frac{1}{Z_n (\beta)}\exp(\beta H_n(S)),
\end{equation}
where 
\begin{equation}
H_n(S)=\sum_{j=1}^n \eta(j,S_j),\textrm{ and } Z_n (\beta)=\textbf{E}\exp(\beta H_n(S)).
\end{equation}
For a given realisation of the environment $\eta$, the measure 
$\g P_{n,\beta}$ gives heavier weight to polymer paths $(i,S_i)_{1\leq i\leq n}$ with lower energy $-H_n(S)$ (configurations of lowest energy are the most probable).
For simplicity we write $\bb E f(\eta)=\bb E f(\eta(0,0))$ for any $f$ such that $f\circ \eta(0,0)$ is integrable. 
Let $\lambda(\beta)=\ln \bb E e^{\beta\eta}$ be the logarithmic moment generating function of $\eta$. We suppose that there exists $B>0$ such that
\begin{equation}\label{hyp}
\bb E e^{B|\eta|}<\infty \textrm{ for } 0\leq \beta\leq B.
\end{equation}
It is well known that the sequence $\bb E\ln Z_n (\beta)$ is
superadditive, hence the limit
\begin{equation}\label{super}
p(\beta)=\Lim n \infty \frac{1}{n}\bb E\ln (Z_n (\beta))=
\sup_n \frac{1}{n}\bb E\ln (Z_n (\beta))\in (-\infty,\lambda(\beta)]
\end{equation}
exists. It is called the free energy of the polymer.
We proved in \cite{LiuWatbled}, \cite{Watbled2012} that if $\bb E e^{\beta|\eta|}<\infty$
for a fixed $\beta>0$, then there exists $K>0$ such that for all $n\geq 1$,
\begin{equation*}
\bb P(\pm\frac{1}{n}(\ln Z_n (\beta)-\bb E \ln Z_n (\beta))>x)
\leq  \left\{\begin{aligned}
&\exp(-\ds\frac{nx^2}{4K}) &\textrm{ if } &x\in (0,2K],\\
&\exp(-n(x-K))  &\textrm{ if }&x\in
(2K, \infty).
\end{aligned}\right.
\end{equation*}
This concentration property implies that under our assumption, $\frac{1}{n}\ln Z_n (\beta)$ converges $\bb P$ a.s. towards $p(\beta)$ for every
$\beta\in (0,B)$. 
This was first proved by Carmona and Hu (\cite{CarmonaHu2002}, Proposition 1.4) for a Gaussian environment,
and by Comets, Shiga and Yoshida (\cite{CometsShigaYoshida}, Proposition 2.5)
for a general environment, but under the condition that $\bb E e^{3\beta|\eta|}<\infty$.\\
We consider the normalized partition function 
$$W_n (\beta)=\frac{Z_n (\beta)}{\bb E Z_n (\beta)}=Z_n (\beta)e^{-n\lambda(\beta)},$$
and the normalized free energy
\begin{equation*}
p_{-}(\beta)=\Lim n \infty \frac{1}{n}\bb E\ln (W_n (\beta))=
p(\beta)-\lambda(\beta)\in (-\infty,0].
\end{equation*}
Bolthausen \cite{Bolthausen} noticed that $(W_n)$ is a positive martingale, hence
converges $\bb P$ a.s. towards a variable $W_\infty$, and that moreover 
$\bb P(W_\infty=0)$ is $0$ or $1$.
When there is no random environment, i.e. when $\beta=0$,  the normalized partition function $W_n$
is constantly equal to one. Accordingly we say that weak disorder holds when $\bb P(W_\infty=0)=0$,
strong disorder holds when $\bb P(W_\infty=0)=1$.
It is immediate that if weak disorder holds, then the normalized free energy 
$p_{-}(\beta)$ equals zero.
Comets and Yoshida proved monotonicity in $\beta$ concerning 
both the dichotomy between weak and strong disorder, and the normalized free energy.
\medskip

\textbf{Theorem} (\cite{CometsYoshida2006}, Theorem 3.2)
\begin{enumerate}
\item 
\textit{There exists a critical value $\beta_0=\beta_c(d)\in [0,\infty]$ 
with $\beta_0=0$ for $d=1,2$, 
$0<\beta_0\leq \infty$ for $d\geq 3$, such that
$\bb P(W_\infty=0)=0$ if $\beta\in \{0\}\cup (0,\beta_0)$,
$\bb P(W_\infty=0)=1$ if $\beta>\beta_0$.}
\item
\textit{The normalized free energy $p_{-}(\beta)$ is non-increasing in $\beta\in [0,\infty)$.
In particular there exists $\beta_c$ with
$\beta_0\leq \beta_c\leq\infty$ such that
$p_{-}(\beta)=0$ if $0\leq \beta\leq \beta_c$,
$p_{-}(\beta)<0$ if $\beta>\beta_c$.}
\end{enumerate}
It is widely believed that the two critical values $\beta_0$ and $\beta_c$ are equal, but it is still an open question in dimension $d$ greater than or equal to three.
The equality $\beta_0=\beta_c=0$ is true in dimension one and two. This was proved by Comets and Vargas (\cite{CometsVargas}) for the dimension one, and by Lacoin (\cite{LacoinNewBounds}) for the dimension two. Moreover in this paper Lacoin improves the result of Comets and Vargas by giving sharp estimates of the normalized free energy at high temperature (that is for small $\beta$).
He proves the following theorem.

\medskip
\textbf{Theorem} (\cite{LacoinNewBounds}, Theorem 1.4 and Theorem 1.5)
\begin{enumerate}
\item 
\textit{When $d=1$ and \eqref{hyp} holds, there exist constants $c$ and $B_0<B$ such that
for every $\beta$ in $[0,B_0)$,}
\begin{equation}\label{genenv}
-\frac{1}{c}\beta^4[1+(\ln\beta)^2]\leq p_{-}(\beta)\leq -c\beta^4.
\end{equation}
\item \textit{When $d=1$ and the environment is Gaussian, then there exists a constant $c$ such that for every $\beta\leq 1$,}
\begin{equation}\label{gauenv}
-\frac{1}{c}\beta^4\leq p_{-}(\beta)\leq -c\beta^4.
\end{equation}
\end{enumerate}
Our aim is to get rid of the logarithmic factor of the lower bound in the case of a general environment.
To prove the lower bound in \eqref{genenv} Lacoin combines the second moment method and a directed percolation argument, whereas in \eqref{gauenv} he uses a specific Gaussian approach similar to what is done in \cite{ToninelliReplica}. 
More precisely his proof relies on two inequalities, both obtained using Gaussian integration by parts. We are able to generalize these inequalities in the case of an infinitely divisible environment.
The first inequality is still obtained by an integration by parts 
formula, valid for infinitely divisible distributions, that we learned in
\cite{CarmonaGuerraHuMejane2006}. Actually, infinite divisibility is not required
for the second inequality, which we prove by using the FKG inequality, in the manner of \cite{CometsYoshida2006}.

\begin{thm}\label{un}
When $d=1$ and $\eta$ has an infinitely divisible distribution satisfying \eqref{hyp},
there exist constants $C>0$ and $B_0<B$ such that
for every $\beta$ in $(0,B_0)$,
\begin{equation}\label{infdivenv}
-C\beta^4\leq p_{-}(\beta).
\end{equation}
\end{thm}

The idea of the proof is to adapt techniques which in the framework of spin glasses are known as replica-coupling and interpolation (introduced in \cite{GuerraToninelli2002}).
Toninelli adapted these techniques to prove disorder irrelevance for the disordered pinning model in the case where the disorder has a Gaussian distribution (\cite{ToninelliReplica}, Theorem 2.1), and mentioned that he can extend his result to the situation where the environment is bounded.
We notice that we can hope to adapt his proof of the Gaussian case to the case where the environment is infinitely divisible, in the same manner as we adapted the proof of Lacoin.
We underline the fact that replica-coupling is not a purely Gaussian technique, and can be adapted to the more general setting of infinite divisible environments.
This was already a key ingredient in \cite{CarmonaGuerraHuMejane2006}. Concerning directed polymers, the Brownian polymer in Poisson environment (\cite{CometsYoshida2005}) might be another good candidate to apply the replica-coupling method.
Let us finally mention that recent work on 1+1 dimensional polymers (for instance \cite{AKQ2012} and \cite{SasamotoSpohn2010}) give another perspective on the $\beta^4$ scaling.

The paper is organized as follows. In the first part we explain the proof
of Theorem \ref{un},
reducing it to the proof of the two inequalities mentioned above.
In the second part we prove the first inequality, using the integration by parts formula for infinitely divisible distributions (Lemma \ref{infdiv}). In the third part we prove the second inequality by using the FKG inequality (Lemma \ref{FKG}).

\section{Proof of Theorem \ref{un}}
We can assume without loss of generality that the variance $\sigma^2$ of $\eta$ is strictly less than one. If it is not the case, we consider $\tilde\eta=\frac{\eta}{2\sigma}$ whose variance is $\frac{1}{4}$ and we use that $p_{-}(\beta)=\tilde p_{-}(2\sigma\beta)$.
We define
\begin{equation}
p_n(\beta)=\frac{1}{n}\bb E\ln W_n(\beta)=\frac{1}{n}\bb E\ln Z_n(\beta)-\lambda(\beta)
\end{equation}
and we recall that 
$$p_{-}(\beta)=\sup_{n\geq 1} p_n(\beta).$$
We use the same strategy as Lacoin in \cite{LacoinNewBounds} and show that
 there exist $B_1>0$ and $c>0$ such that for every $n\geq 1$ and every $\beta\in [0,B_1)$,
 \begin{equation}\label{Wat}
 p_n(\beta)\geq (1-e^c)\frac{1}{2n}\ln E^{\otimes 2}e^{2\beta^2L_n(S^1,S^2)},
 \end{equation}
  where $S^1$, $S^2$ are two independent copies of the simple random walk, $E^{\otimes 2}$ is the expectation on the product space $(\Sigma^2,\cal E^{\otimes 2})$,  
  and 
$$L_n(S^1,S^2)=\sum_{i=1}^n \mathbf 1_{S_i^1=S_i^2}$$
is the number of intersections of the two directed paths $(i,S_i^1)_{1\leq i\leq n}$
and $(i,S_i^2)_{1\leq i\leq n}$.
Now $E^{\otimes 2}e^{2\beta^2L_n(S^1,S^2)}$ is the partition function of a homogeneous pinning system of size $n$ and parameter $2\beta^2$ with underlying renewal process the sets of zeros of the random walk $S^1-S^2$.
We refer to \cite{Giacomin2007}, section 1.2
 for a survey of the homogeneous pinning model, where it is proved (Proposition 1.1 and Remark 1.2) that 
 $$\lim_{n\to\infty} \frac{1}{n}\ln E^{\otimes 2}e^{tL_n(S^1,S^2)}=F(t),
 \textrm{ with }F(t)\sim \frac{t^2}{2}\textrm{ when }t\to 0^+.$$
 Letting $n$ tend to infinity in \eqref{Wat} we find that
 $$p_{-}(\beta)\geq \frac{1-e^c}{2}F(2\beta^2).$$
 This implies that for any $C>e^c-1$, there exists $B_0<B$ such that
for all $0\leq \beta< B_0$,
$$p_{-}(\beta)\geq -C\beta^4,$$
which proves the lower bound result \eqref{infdivenv}.
To prove \eqref{Wat} we proceed as in \cite{LacoinNewBounds} by interpolation between the two functions
$$p_n(\beta)=\frac{1}{n}\bb E \ln Ee^{\beta H_n(S)-n\lambda(\beta)}
 =\frac{1}{2n}\bb E \ln E^{\otimes 2}e^{\beta H_n(S^1,S^2)-2n\lambda(\beta)},$$
 where $$H_n(S^1,S^2)=H_n(S^1)+H_n(S^2),$$ and
 $$F_n(\beta)=\frac{1}{2n}\ln E^{\otimes 2}e^{2\beta^2L_n(S^1,S^2)}
 =\frac{1}{2n}\bb E\ln E^{\otimes 2}e^{2\beta^2L_n(S^1,S^2)}.$$
Let $\beta$ be fixed in $]0,B[$.
We define, for $t\in [0,1]$ and $u\geq 0$,
$$\varphi_n(t,u)=
\frac{1}{2n}\bb E\ln E^{\otimes 2} e^{\sqrt t \beta H_n(S^1,S^2)-2n\lambda(\sqrt t \beta)+u\beta^2L_n(S^ 1,S^ 2)},$$
$$\phi_n(t)=\varphi_n(t,0)=\frac{1}{n}\bb E\ln E e^{\sqrt t \beta H_n(S)-
n\lambda(\sqrt t \beta)}=p_n(\sqrt t\beta),$$
so that
$$\phi_n(1)=\varphi_n(1,0)=p_n(\beta),\ \  \varphi_n(0,2)=F_n(\beta).$$
The inequality \eqref{Wat} amounts to
\begin{equation}\label{Wat2}
\phi_n(1)\geq (1-e^c)\varphi_n(0,2).
\end{equation}
We shall prove that 
\begin{equation}\label{Gron}\forall t\in [0,1], \ \ \phi'_n(t)\geq c(\phi_n(t)-\varphi_n(0,2)),
\end{equation}
which by Gronwall's lemma implies \eqref{Wat2}.
We write
\begin{equation}\label{sum}
\phi_n(t)-\varphi_n(0,2)
=\varphi_n(t,0)-\varphi_n(t,2-t)+\varphi_n(t,2-t)-\varphi_n(0,2)
\end{equation}
and we consider separately the two differences of the right member. 
To prove \eqref{Gron} it is enough to prove the two following inequalities.
\begin{equation}\label{first}
\exists c>0,\quad\forall t\in [0,1], \ \ \phi'_n(t)\geq c(\varphi_n(t,0)-\varphi_n(t,2-t)),
\end{equation}
\begin{equation}\label{second}\forall t\in [0,1], \ \ \
\varphi_n(t,2-t)-\varphi_n(0,2)\leq 0.
\end{equation}
To prove \eqref{first} we need infinite divisibility whereas \eqref{second} is true 
in general.
Let us first prove \eqref{first}. 
The function $u\mapsto \varphi_n(t,u)$ is convex non decreasing, hence for every $t\in[0,1]$,
\begin{equation}\label{sum1}\varphi_n(t,2-t)-\varphi_n(t,0)\geq (2-t)\frac{\partial\varphi_n}{\partial u}(t,0)\geq \frac{\partial\varphi_n}{\partial u}(t,0).
\end{equation}
A simple calculation shows that
\begin{equation}\label{cal1}\frac{\partial\varphi_n}{\partial u}(t,0)=\frac{\beta^2}{2n}\bb E E_{n,\sqrt t\beta}^{\otimes 2}L_n(S^1,S^2),\end{equation}
where $E_{n,\beta}$
is the expectation with respect to the probability measure $P_{n,\beta}$,
and $E_{n,\beta}^{\otimes 2}$ is the expectation with respect to the product measure 
$P_{n,\beta}^{\otimes 2}$ on the product space $(\Sigma^2,\cal E^{\otimes 2})$.
On the other hand,
\begin{equation}\label{cal2}\phi_n'(t)=\frac{\beta}{2\sqrt t}p_n'(\sqrt t\beta).
\end{equation}
Now the quantities \eqref{cal1} and \eqref{cal2} are related through the following lemma, whose proof we postpone to the next section.
\begin{lem}\label{infdiv}
If the law of $\eta$ is infinitely divisible, then there exists $c>0$ such that for every $\beta\in [0,\frac{B}{2}),$
\begin{equation}\label{carmona}
p_n'(\beta)\geq -\frac{c\beta}{n}\bb E E_{n,\beta}^{\otimes 2}L_n(S^1,S^2).
\end{equation}
\end{lem}
We use successively \eqref{cal2}, the Lemma \ref{infdiv}, \eqref{cal1}, and \eqref{sum1},
and we obtain that provided we take $\beta<\frac{B}{2}$,
there exists $c>0$ such that for every $t\in[0,1]$,
\begin{equation*}
\phi_n'(t)\geq -\frac{c\beta^2}{2n}\bb E E_{n,\sqrt t\beta}^{\otimes 2}L_n(S^1,S^2)
\geq c(\varphi_n(t,0)-\varphi_n(t,2-t)),\end{equation*}
 which is exactly \eqref{first}.
Now we prove the second inequality \eqref{second} thanks to the following lemma.
\begin{lem}\label{FKG}
If $\sigma^2<1$ then there exists $B_2>0$ such that for every fixed $\beta$ in $(0,B_2)$, 
\begin{equation}\label{point2bis}
\frac{\partial}{\partial u}\varphi_n(t,u)-\frac{\partial}{\partial t}\varphi_n(t,u)\geq 0 \ \ \ \textrm{ for } t\in [0,1] \textrm{ and }u\geq 0.
\end{equation}
\end{lem}
We postpone the proof of Lemma \ref{FKG} and we conclude the proof of the theorem.
We set $B_1=\min(\frac{B}{2},B_2)$ and fix $\beta$ in $[0,B_1)$.
Then according to Lemma \ref{FKG}
\begin{equation*}
\varphi_n(t,2-t)-\varphi_n(0,2)=\int_0^t
\frac{\partial\varphi_n}{\partial t}(s,2-s)-\frac{\partial\varphi_n}{\partial u}(s,2-s)ds\leq 0,
\end{equation*}
which is inequality \eqref{second}.
Combining \eqref{sum}, \eqref{first}, and \eqref{second}, we obtain \eqref{Gron}, which ends the proof.

\begin{rem}
If the law of $\eta$ is standard normal, then the inequality \eqref{carmona}
is an equality with $c=1$, and the inequality \eqref{point2bis} holds true.
Both facts are showned in \cite{LacoinNewBounds} by using Gaussian integration by parts.
\end{rem}

\section{Proof of Lemma \ref{infdiv}}

According to \cite{CarmonaGuerraHuMejane2006}, Proposition 11, there exists $c_2>0$ depending on $\beta$ such that
\begin{equation}\label{carmonac2}
p_n'(\beta)\geq -\frac{c_2}{n}\bb E E_{n,\beta}^{\otimes 2}L_n(S^1,S^2).
\end{equation}
We recall rapidly the proof of \cite{CarmonaGuerraHuMejane2006}.
We write
$$
np_n'(\beta)=\bb E\frac{EH_n(S)e^{\beta H_n(S)}}{E e^{\beta H_n(S)}}-n\lambda'(\beta)
=\sum_{(i,x)} \bb E\eta(i,x)f_{i,x}(\eta(i,x))-n\lambda'(\beta).
$$
The distribution of $\eta$ is infinitely divisible, therefore we have a Lévy Khinchine formula
$$\lambda(\beta)=c_0\beta+\frac{\sigma^2}{2}\beta^2
+\int(e^{\beta u}-1-\beta u\mathbf 1_{\abs{u}\leq 1})\pi(du),$$
where $c_0\in\bb R$, $\sigma\geq 0$, and $\pi$ is a measure on $\bb R^*$ such that $1\wedge u^2\in L^1(\pi)$.
With these notations we have the following lemma.
\medskip

\textbf{Lemma} (\cite{CarmonaGuerraHuMejane2006}, Lemma 10).
\textit{For any bounded differentiable $f$ with bounded derivative, one has the following integration by parts formula:}
$$\bb E\eta f(\eta)=c_0\bb E f(\eta)+\sigma^2 \bb E f'(\eta)+\int_\bb R (\bb Ef(\eta +u)- \mathbf 1_{\abs{u}\leq 1}\bb Ef(\eta))u\pi(du).$$
Applying the lemma, Carmona and al. calculate
$$np_n'(\beta)= -\sigma^2\beta\bb E E_{n,\beta}^{\otimes 2}L_n(S^1,S^2)
-\sum_{(i,x)}\bb E(P_{n,\beta}(S_i=x))^2\int_\bb R 
\frac{(e^{\beta u}-1)e^{\beta u}u}{(e^{\beta u}-1)\mu_n(S_i=x)+1}\pi(du).$$
Noting that
$$\sum_{(i,x)}(P_{n,\beta}(S_i=x))^2=E_{n,\beta}^{\otimes 2}L_n(S^1,S^2),$$
and that
\begin{multline}\label{carmona2}
\sup_{0\leq a\leq 1} 
\int_\bb R \frac{(e^{\beta u}-1)e^{\beta u}u}{(e^{\beta u}-1)a+1}\pi(du)\\
\leq \int_{-\infty}^0 \abs{u}(1-e^{\beta u})\pi(du)+
\int_{0}^{+\infty} ue^{\beta u}(e^{\beta u}-1)\pi(du)=c'(\beta),
\end{multline}
which is finite provided that $\bb E e^{2\beta\eta}$ is finite,
they conclude that \eqref{carmonac2} is true with $c_2=c'(\beta)+\sigma^2\beta$.

Now for $u<0$ we have $1-e^{\beta u}\leq\beta\abs{u}$,
and for $u>0$ we have
$e^{\beta u}-1\leq \beta ue^{\beta u},$
so returning to \eqref{carmona2} we deduce that
for every $\beta\in[0,\frac{B}{2}]$, 
$$c'(\beta)
\leq \beta\left(\int_{-\infty}^0 u^2\pi(du)+\int_{0}^{+\infty} 
u^2e^ {B u}\pi(du)\right),$$
and conclude that \eqref{carmona} is satisfied with
$$c=\sigma^ 2+\int_{-\infty}^0 u^2\pi(du)+\int_{0}^{+\infty} u^2e^ {B u}\pi(du).$$

\section{Proof of Lemma \ref{FKG}}

We calculate
\begin{equation}\label{F1}
\frac{\partial}{\partial u}\varphi_n(t,u)-\frac{\partial}{\partial t}\varphi_n(t,u)
=\frac{1}{2n}E^{\otimes 2}[e^{u\beta^2L_n(S^ 1,S^ 2)}\bb I(S^1,S^2)],
\end{equation}
where
$$
\begin{aligned}
\bb I(S^1,S^2)&=\bb E\frac{ e^{\sqrt t \beta H_n(S^1,S^2)-2n\lambda(\sqrt t \beta)}
(\beta^2L_n(S^ 1,S^ 2)-
\frac{\beta}{2\sqrt t}H_n(S^1,S^2)+\frac{\beta n}{\sqrt t}\lambda'(\sqrt t\beta))}
{E^{\otimes 2} e^{\sqrt t \beta H_n(S^1,S^2)-2n\lambda(\sqrt t \beta)+u\beta^2L_n(S^ 1,S^ 2)}}\\
&=\bb E\frac{ e^{\sqrt t \beta H_n(S^1,S^2)}
(\beta^2L_n(S^ 1,S^ 2)-
\frac{\beta}{2\sqrt t}H_n(S^1,S^2)+\frac{\beta n}{\sqrt t}\lambda'(\sqrt t\beta))}
{E^{\otimes 2} e^{\sqrt t \beta H_n(S^1,S^2)+u\beta^2L_n(S^ 1,S^ 2)}}
\end{aligned}$$
We show that $\bb I(S^1,S^2)$ is non negative for every $(S^1,S^2)$ in $\Sigma^2$.
Let $(S^1,S^2)$ be fixed in $\Sigma^2$. 
We define
$$I=\{i\in\interventier{1}{n}; S^1_i=S^2_i\}\textrm{ and } 
A=\{(i,x);1\leq i\leq n\textrm{ and } x=S^1_i\textrm{ or }S^2_i\}.$$
We write $H_n$ for $H_n(S^1,S^2)$ and $L_n$ for $L_n(S^1,S^2)$. For 
$v=(v_{(i,x)})_{(i,x)\in A}$, we define
$$H_n(v)=\sum_{i=1}^{n}v_{(i,S^1_i)}+v_{(i,S^2_i)}=
\sum_{i\in I}2v_{(i,S^1_i)}+\sum_{i\notin I} v_{(i,S^1_i)}+v_{(i,S^2_i)} $$
Then the measure defined by
$$\mu(dv)=\frac{e^{\sqrt t \beta H_n(v)}}{\bb E(e^{\sqrt t \beta H_n})} \bb P_\eta^{\otimes(2n-L_n)}$$
is a product probability measure on $\bb R^{2n-L_n}$, and 
we can write $\bb I(S^1,S^2)$ as an expectation with respect to $\mu$:
\begin{multline}\label{F2}
\bb I(S^1,S^2)
=\int \frac{\beta^2L_n-
\frac{\beta}{2\sqrt t}H_n(v)+\frac{\beta n}{\sqrt t}\lambda'(\sqrt t\beta)}
{E^{\otimes 2} e^{\sqrt t \beta H_n(v)+u\beta^2L_n}}\bb E(e^{\sqrt t \beta H_n})\mu(dv)\\
=\int X(v)Y(v)\mu(dv),
\end{multline}
with
$$X(v)=\frac{1}
{E^{\otimes 2} e^{\sqrt t \beta H_n(v)+u\beta^2L_n}}$$
and
$$Y(v)=\Big(\beta^2 L_n-
\frac{\beta}{2\sqrt t}H_n(v)+\frac{\beta n}{\sqrt t}\lambda'(\sqrt t\beta)\Big)\bb E(e^{\sqrt t \beta H_n}).$$
As $X$ and $Y$ are decreasing random variables in $L^2(\mu)$, we deduce from the FKG inequality that
\begin{equation}\label{F3}
\int X(v)Y(v)\mu(dv)\geq \int X(v)\mu(dv)\times\int Y(v)\mu(dv).
\end{equation}
We refer to \cite{Grimmett1999}, Theorem 2.4 or \cite{Liggett1985}, chapter II for 
more information about the FKG inequality.
We calculate $\int Y(v)\mu(dv)$ and find that
$$\int Y(v)\mu(dv)
=\Big(\beta^2L_n+\frac{\beta n}{\sqrt t}\lambda'(\sqrt t\beta)\Big)
\bb E (e^{\sqrt t \beta H_n})
-\frac{\beta}{2\sqrt t}\bb E(H_n e^{\sqrt t \beta H_n}).$$
Using that $\bb E(\eta e^{u\eta})=\lambda'(u)e^{\lambda(u)}$, we calculate
$$\bb E(H_n e^{\sqrt t \beta H_n})=
\Big(2L_n \lambda'(2\sqrt t\beta)+2(n-L_n)\lambda'(\sqrt t\beta)\Big)\bb E(e^{\sqrt t \beta H_n}),
$$
hence
$$\int Y(v)\mu(dv)
=\beta^2L_n\Big(1-\frac{\lambda'(2\sqrt t\beta)-\lambda'(\sqrt t\beta)}{\sqrt t\beta}\Big)
\bb E e^{\sqrt t \beta H_n}.$$
There exists $a\in (\sqrt t\beta,2\sqrt t\beta)$ such that 
$\frac{\lambda'(2\sqrt t\beta)-\lambda'(\sqrt t\beta)}{\sqrt t\beta}=\lambda''(a)$.
Since $\lambda''(0)=\sigma^2$ is strictly less than one there exists $B_2>0$ such that 
for every $s$ in $[0,2B_2)$, $\lambda''(s)$ remains strictly less than one.
If $\beta\in (0,B_2)$ then $a\in (0,2B_2)$, hence $\lambda''(a)< 1$, from which we deduce that $\int Y(v)\mu(dv)$ is non negative.
As $\int X(v)\mu(dv)$ is also obviously non negative, recalling \eqref{F1}, \eqref{F2}
and \eqref{F3}, we get the result.

\section*{Acknowledgements}
The author thanks very much the anonymous referee for his careful reading and valuable comments. 

\providecommand{\bysame}{\leavevmode\hbox to3em{\hrulefill}\thinspace}
\providecommand{\MR}{\relax\ifhmode\unskip\space\fi MR }
\providecommand{\MRhref}[2]{%
  \href{http://www.ams.org/mathscinet-getitem?mr=#1}{#2}
}
\providecommand{\href}[2]{#2}

\end{document}